\documentclass{elsart}

\usepackage{graphicx}
\usepackage{amssymb}
\usepackage{psfrag}
\usepackage{here}

\newcommand{\real}{\mathbb{R}}
\newcommand{\cgq}{$\mathrm{cG}(q)$}
\newcommand{\dgq}{$\mathrm{dG}(q)$}
\newcommand{\mcgq}{$\mathrm{mcG}(q)$}
\newcommand{\mdgq}{$\mathrm{mdG}(q)$}
\newcommand{\PD}[2]{\frac{\partial #1}{\partial #2}}
\newcommand{\diag}{\mathrm{diag}}

\begin{document}

\begin{frontmatter}

\title{Multi-Adaptive Time Integration}
\author{Anders Logg}
\ead{logg@math.chalmers.se}
\ead[url]{http://www.math.chalmers.se/\~{}logg}
\address{Department of Computational Mathematics\\Chalmers University of Technology\\ G\"oteborg, Sweden}

\begin{abstract}
Time integration of ODEs or time-dependent PDEs with required resolution of
the fastest time scales of the system, can be very costly if
the system exhibits multiple time scales of different magnitudes. If the
different time scales are localised to different components,
corresponding
to localisation in space for a PDE, efficient time integration thus
requires that we use different time steps for different components.

We present an overview of the multi-adaptive Galerkin methods \mcgq\ and
\mdgq\ recently introduced in a series of papers by the author. In these
methods,
the time step sequence is selected individually and adaptively for each
component, based on an a posteriori error estimate of the global error.

The multi-adaptive methods require the solution of large systems of
nonlinear algebraic equations which are solved using explicit-type
iterative solvers (fixed point iteration).
If the system is stiff, these iterations may fail to
converge, corresponding to the well-known fact
that
standard explicit methods are inefficient for stiff systems. To resolve
this problem, we present an adaptive strategy for explicit time
integration of stiff ODEs, in which the explicit method is adaptively
stabilised by a small number of small, stabilising time steps.
\end{abstract}

\end{frontmatter}

\section{Introduction}

In earlier work \cite{logg:article:ma1,logg:article:ma2}, we have introduced the
multi-adaptive Galerkin methods \mcgq\ and \mdgq\
for ODEs of the type
\begin{equation}
  \left\{
    \begin{array}{rcl}
      \dot{u}(t) &=& f(u(t),t), \quad t\in(0,T], \\
      u(0) &=& u_0,
    \end{array}
  \right.
  \label{eq:u'=f}
\end{equation}
where $u : [0,T] \rightarrow \real^N$ is the solution to be computed,
$u_0 \in \real^N$ a given initial condition,
$T>0$ a given final time,
and $f : \real^N \times (0,T] \rightarrow \real^N$ a given
function that is Lipschitz-continuous in $u$ and bounded.
We use the term \emph{multi-adaptivity} to describe methods with individual time-stepping
for the different components $u_i(t)$ of the solution vector $u(t)=(u_i(t))$, including
(i) time step length,
(ii) order, and
(iii)
quadrature, all chosen adaptively in a computational feed-back process.

Surprisingly, individual time-stepping for ODEs has received little attention in the
large literature on numerical methods for ODEs, see e.g.
\cite{dahlquist:thesis,hairerwanner:book1,hairerwanner:book2,butcher:book,shampine:book}.
For specific applications, such as the $n$-body problem, methods with
individual time-stepping have been used, see e.g.
\cite{makino:local,alexander:local,dave:local},
but a general methodology has been
lacking.
For time-dependent PDEs, in particular for conservation laws of the type
$\dot{u} + f(u)_x = 0$, attempts have been made to construct methods with individual
(locally varying in space) time steps. Flaherty et al. \cite{FlaLoy97} have constructed a method
based on the discontinuous Galerkin method combined with local explicit Euler time-stepping.
A similar approach is taken in \cite{DawKir01} where a method based on the original work by
Osher and Sanders \cite{OshSan83} is
presented for conservation laws in one and two space dimensions. Typically the time steps used
are based on local CFL conditions rather than error estimates for the global error and
the methods are low order in time (meaning $\leq 2$).
We believe that our work on multi-adaptive Galerkin methods (including error estimation and arbitrary order methods)
presents a general methodology to
individual time-stepping, which will result in efficient integrators both for ODEs and
time-dependent PDEs.

The multi-adaptive methods are developed within the general framework of adaptive
Galerkin methods based on piecewise polynomial approximation (finite element methods)
for differential equations, including the continuous and discontinuous Galerkin methods \cgq\ and \dgq\,
which we extend to their multi-adaptive analogues \mcgq\ and \mdgq.
Earlier work on adaptive error control for the \cgq\ and \dgq\ methods include
\cite{delfour:dg,ejt:dg,claes:dg,estep:cg,estep:dg,EstWil96}.
See also \cite{EJ:parab:I,EJ:parab:II,EJ:parab:III,EJ:parab:IV,EJ:parab:V,EJ:parab:VI},
and \cite{EJ:actanumerica} or \cite{EstLar00} in particular for an overview of adaptive error control based on duality techniques.
The approach to error analysis and adaptivity presented in these references naturally carries over to the multi-adaptive methods.

\subsection{The stiffness problem}

The classical wisdom developed in the 1950s regarding stiff ODEs
is that efficient
integration requires implicit (A-stable) methods, at
least outside transients where the time steps may be chosen
large from accuracy point of view. Using an explicit method
(with a bounded stability region) the time steps have to
be small at all times for stability reasons, in particular outside
transients, and the advantage of a low cost per time step for the
explicit method is counter-balanced by the necessity of taking a
large number of small time steps. As a result, the overall efficiency
of an explicit method for a stiff ODE is small.

We encounter the same problem when we try to use explicit fixed point
iteration to solve the discrete equations given by the multi-adaptive
Galerkin methods \mcgq\ and \mdgq. However, it turns out that if a
sequence of large (unstable) time steps are accompanied by a suitable
(small) number of small time steps, a stiff system can be stabilised
to allow integration with an effective time step much larger than the
largest stable time step given by classical stability analysis. This
idea of stabilising a stiff system using the inherent damping property
of the stiff system itself was first developed in an automatic and adaptive setting
in \cite{logg:article:stiffode}, and will be further explored
in the full multi-adaptive setting. A
similar approach is taken in recent independent work by Gear and
Kevrekidis \cite{GeaKev02}. The relation to Runge-Kutta methods based
on Chebyshev polynomials discussed by Verwer in \cite{Ver96} should
also be noted.

\subsection{Notation}

\begin{figure}[htbp]
        \begin{center}
                \psfrag{0}{$0$}
                \psfrag{i}{$i$}
                \psfrag{k}{$k_{ij}$}
                \psfrag{K}{$K_n$}
                \psfrag{T}{$T$}
                \psfrag{I}{$I_{ij}$}
                \psfrag{t1}{\small $t_{i,j-1}$}
                \psfrag{t2}{\small $t_{ij}$}
                \psfrag{T1}{$T_{n-1}$}
                \psfrag{T2}{$T_n$}
                \includegraphics[width=12cm]{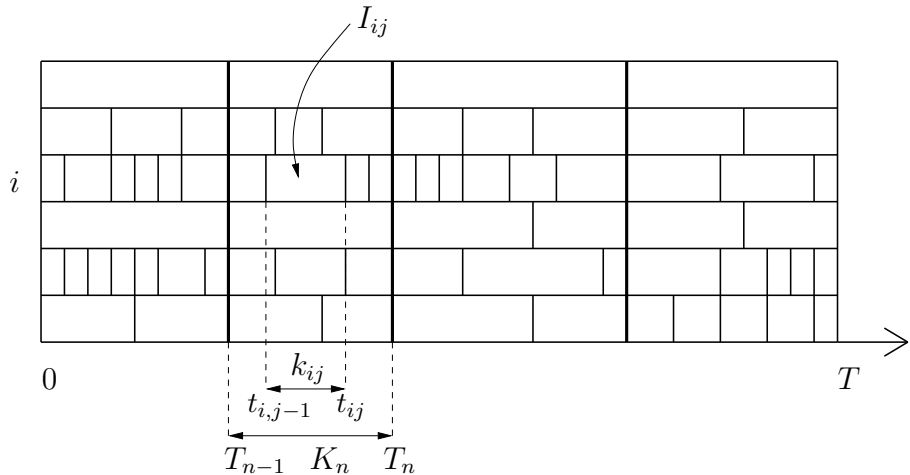}

                \caption{Individual partitions of the interval \ $(0,T]$ for different components. Elements
                                        between common synchronised time levels are organised
                                        in time slabs. In this example, we have $N=6$ and $M=4$.}
                \label{fig:notation}
        \end{center}
\end{figure}

The following notation is used in the discussion of the multi-adaptive
Galerkin methods below:
Each component $U_i(t)$, $i=1,\ldots,N$, of the approximate
$\mathrm{m(c/d)G}(q)$ solution $U(t)$ of (\ref{eq:u'=f}) is a
piecewise polynomial on a partition of $(0,T]$ into $M_i$ subintervals.
Subinterval $j$ for component $i$ is denoted by $I_{ij}=(t_{i,j-1},t_{ij}]$,
and the length of the subinterval is given by the local \emph{time step} $k_{ij}=t_{ij}-t_{i,j-1}$.
This is illustrated in Figure \ref{fig:notation}.
On each subinterval $I_{ij}$, $U_{i}\vert_{I_{ij}}$ is a polynomial
of degree $q_{ij}$ and we refer to $(I_{ij},U_i\vert_{I_{ij}})$ as an \emph{element}.

Furthermore, we shall assume that the interval $(0,T]$ is
partitioned into blocks between certain synchronised time levels
$0=T_0<T_1<\ldots<T_M=T$. We refer to the set of intervals $\mathcal{T}_n$ between two synchronised time levels $T_{n-1}$ and $T_n$ as a \emph{time slab}:
$\mathcal{T}_n = \{ I_{ij} : T_{n-1} \leq t_{i,j-1} < t_{ij} \leq T_n \}$, and we denote
the length of a time slab by $K_n = T_n - T_{n-1}$.
The partition consisting of the entire collection of intervals is denoted by $\mathcal{T}=\cup \mathcal{T}_n$.

\subsection{Outline}

The outline of the paper is as follows: In Section \ref{sec:methods},
we formulate the multi-adaptive Galerkin methods \mcgq\ and \mdgq. In
Section \ref{sec:errorcontrol}, we discuss error control and
adaptivity. In particular, we show how
to choose the individual time steps based on an a posteriori error
estimate for the global error. In Section \ref{sec:iterative}, we give
a quick overview of an iterative method (based on fixed point
iteration) for the system of nonlinear discrete equations that needs
to be solved on each time slab, and in Section \ref{sec:stiff} we describe a technique
that can be used to stabilise the explicit fixed point iterations for
stiff problems. Finally, in Section \ref{sec:examples}, we present
a number of numerical examples chosen to illustrate both the potential
of multi-adaptivity and the use of explicit fixed point iteration
(or explicit time-stepping) for stiff problems.

\section{Multi-Adaptive Galerkin}
\label{sec:methods}

\subsection{Multi-adaptive continuous Galerkin, \mcgq}

To give the definition of the \mcgq\ method, we define the
\emph{trial space} $V$ and the \emph{test space} $W$ as
\begin{equation}
  \begin{array}{rcl}
    V &=& \{v \in [\mathcal{C}([0,T])]^N : v_i|_{I_{ij}}\in \mathcal{P}^{q_{ij}}(I_{ij}), \     j=1,\ldots,M_i, \ i=1,\ldots,N \},\\
    W &=& \{v : v_i|_{I_{ij}}\in \mathcal{P}^{q_{ij}-1}(I_{ij}), \
    j=1,\ldots,M_i, \ i=1,\ldots,N \},\\
  \end{array}
  \label{eq:spaces,mcg}
\end{equation}
where $\mathcal{P}^q(I)$ denotes the linear space of polynomials of
degree $q\geq 0$ on the interval $I$. In other words, $V$ is the space
of continuous piecewise polynomials of
degree $q=q_i(t)=q_{ij}, \ t\in I_{ij}$ on the partition $\mathcal{T}$,
and $W$ is the space of (in general discontinuous) piecewise polynomials
of degree $q-1$ on the same partition.

We define the \mcgq\ method for (\ref{eq:u'=f}) as follows:
Find $U\in V$ with $U(0)=u_0$, such that
\begin{equation}
         \int_0^T (\dot{U},v)\ dt = \int_0^T (f(U,\cdot),v) \ dt \quad \forall v\in W,
    \label{eq:fem,mcg}
\end{equation}
where $(\cdot,\cdot)$ denotes the standard inner product in $\real^N$.
If now for each local interval $I_{ij}$ we take $v_n=0$ when $n\neq i$ and
$v_i(t)=0$ when $t\not\in I_{ij}$, we can rewrite the global problem
(\ref{eq:fem,mcg}) as a number of successive local problems for each
component: For $i=1,\ldots,N$, $j=1,\ldots,M_i$,
find $U_i|_{I_{ij}}\in \mathcal{P}^{q_{ij}}(I_{ij})$ with $U_i(t_{i,j-1})$ given from the previous time interval, such that
\begin{equation}
    \int_{I_{ij}} \dot{U}_i v \ dt = \int_{I_{ij}} f_i(U,\cdot) v \ dt \quad \forall v\in \mathcal{P}^{q_{ij}-1}(I_{ij}).
    \label{eq:fem,mcg,local}
\end{equation}
We define the \emph{residual} $R$ of the approximate solution $U$
to be $R(U,t) = \dot{U}(t) - f(U(t),t)$. In terms of the residual, we can rewrite (\ref{eq:fem,mcg,local}) as
$\int_{I_{ij}} R_i(U,\cdot) v \ dt = 0$ for all $v\in \mathcal{P}^{q_{ij}-1}(I_{ij})$,
i.e., the residual is orthogonal to the test space on every local interval. We refer
to this as the \emph{Galerkin orthogonality} of the \mcgq\ method.

\subsection{Multi-adaptive discontinuous Galerkin, \mdgq}

For the \mdgq\ method, we define the trial and test spaces by
\begin{equation}
    V = W = \{v : v_i|_{I_{ij}}\in \mathcal{P}^{q_{ij}}(I_{ij}),
    j=1,\ldots,M_i, \ i=1,\ldots,N \},
  \label{eq:spaces,mdg}
\end{equation}
i.e., both trial and test functions are (in general discontinuous)
piecewise polynomials of degree $q=q_i(t)=q_{ij}, \ t\in I_{ij}$ on
the partition $\mathcal{T}$.

We define the \mdgq\ method for (\ref{eq:u'=f}) as follows,
similar to the definition of the continuous method:
Find $U\in V$ with $U(0^-)=u_0$, such that
\begin{equation}
  \sum_{i=1}^N \sum_{j=1}^{M_i}
  \left[
    [U_i]_{i,j-1} v(t_{i,j-1}^+) + \int_{I_{ij}} \dot{U}_i v_i \ dt
  \right] =
  \int_0^T (f(U,\cdot),v) \ dt \quad \forall v\in W.
  \label{eq:dual,dg}
\end{equation}
where $[\cdot]$ denotes the jump,
i.e., $[v]_{ij} = v(t_{ij}^+) - v(t_{ij}^-)$.

The \mdgq\ method in local form, corresponding to
(\ref{eq:fem,mcg,local}), reads:
For $i=1,\ldots,N$, $j=1,\ldots,M_i$, find $U_i|_{I_{ij}}\in \mathcal{P}^{q_{ij}}(I_{ij})$, such that
\begin{equation}
    [U_i]_{i,j-1} v(t_{i,j-1}) + \int_{I_{ij}} \dot{U}_i v \ dt = \int_{I_{ij}} f_i(U,\cdot) v \ dt \quad \forall v\in \mathcal{P}^{q_{ij}}(I_{ij}),
    \label{eq:fem,mdg,local}
\end{equation}
where the initial condition is specified
for $i=1,\ldots,N$, by $U_i(0^-) = u_i(0)$.
In the same way as for the continuous method, we define the residual $R$ of the approximate
solution $U$ to be $R(U,t)=\dot{U}(t) - f(U(t),t)$, defined on the inner of every local
interval $I_{ij}$, and rewrite (\ref{eq:fem,mdg,local}) in the form
$[U_i]_{i,j-1} v(t_{i,j-1}) + \int_{I_{ij}} R_i(U,\cdot) v \ dt = 0$ for all $v\in \mathcal{P}^{q_{ij}}(I_{ij})$.
We refer to this as the Galerkin orthogonality of the \mdgq\ method.

\section{Error Control and Adaptivity}
\label{sec:errorcontrol}

Our goal is to compute an approximation $U(T)$ of the exact solution
$u(T)$ of (\ref{eq:u'=f}) at final time $T$ within a given tolerance $\mathrm{TOL}>0$,
using a minimal amount of computational work. This goal includes an
aspect of \emph{reliability} (the error should be less than the
tolerance) and an aspect of \emph{efficiency} (minimal computational
work). To measure the error we choose a norm, such as the Euclidean
norm $\|\cdot\|$ on $\real^N$, or more generally some other quantity
of interest.

We discuss below both a priori and a posteriori error estimates
for the multi-adaptive Galerkin methods, and the application of the
a posteriori error estimates in multi-adaptive time-stepping.

\subsection{A priori error estimates}

Standard (duality-based) a priori error estimates show that the order for
the ordinary Galerkin methods $\mathrm{cG}(q)$ and $\mathrm{dG}(q)$ is
$2q$ and $2q+1$, respectively. A generalisation of these
estimates to the multi-adaptive methods gives the same result. The
multi-adaptive continuous Galerkin method $\mathrm{mcG}(q)$ is thus of
order $2q$, and the multi-adaptive discontinuous Galerkin method
$\mathrm{mdG}(q)$ is of order $2q+1$.

\subsection{A posteriori error estimates}

A posteriori error analysis in the general framework of \cite{EJ:actanumerica}
relies on the concept of the \emph{dual problem}. The dual problem
of the initial value problem (\ref{eq:u'=f}) is the linearised backward
problem given by
\begin{equation}
  \left\{
    \begin{array}{rcl}
      -\dot{\phi} &=& J^*(u,U,\cdot)\phi \quad \mbox{ on } [0,T),\\
      \phi(T) &=& e(T)/\|e(T)\|,
    \end{array}
  \right.
  \label{eq:dual,continuous}
\end{equation}
where the Jacobian $J$ is given by
$J(u,U,\cdot{}) = \int_0^1 \PD{f}{u}(s u + (1-s) U,\cdot{}) \ ds$
and $^*$ denotes the transpose. We use the dual problem to represent the error in terms
of the dual solution $\phi$ and the residual $R$. For the \mcgq\ method
the representation formula is given by
\begin{equation}
	\| e(T) \| = \int_0^T (R,\phi) \ dt,
	\label{eq:representation,cg}
\end{equation}
and for the \mdgq\ method, we obtain
\begin{equation}
	\| e(T) \| = \sum_{i=1}^N \sum_{j=1}^{M_i} \ [U_i]_{i,j-1} \phi_i(t_{i,j-1}) + \int_{I_{ij}} R_i(U,\cdot) \phi_i \ dt.
	\label{eq:representation,dg}
\end{equation}

Using the Galerkin orthogonalities together with special interpolation estimates (see \cite{logg:article:ma1}),
we obtain a posteriori error estimates of the form
\begin{equation} \label{eq:estimate,cg}
	\| e(T) \| \leq \sum_{i=1}^N S^{[q_i]}_i \max_{[0,T]} \left\{ C k_{i}^{q_{i}} r_{i} \right\},
\end{equation}
for the \mcgq\ method, and
\begin{equation} \label{eq:estimate,dg}
	\| e(T) \| \leq \sum_{i=1}^N S^{[q_i+1]}_i \max_{[0,T]} \left\{ C k_{i}^{q_{i}+1} r_{i} \right\},
\end{equation}
for the \mdgq\ method, where $C$ is an interpolation constant, $r_i$
is a local measure of the residual, and the individual \emph{stability
factors} $S_i$ are given by $S_i^{[q_i]} = \int_0^T \vert
\phi_i^{(q_i)} \vert \ dt$.  Typically, the stability factors are
of moderate size for a stiff problem (and of unit size for a parabolic
problem), which means that accurate computation is possible over long
time intervals. Note that the Lipschitz constant, which is large
for a stiff problem, is not present in these estimates.

The analysis can be extended to include also \emph{computational errors}, arising from solving
the discrete equations using an iterative method, and \emph{quadrature errors},
arising from evaluating the integrals in (\ref{eq:fem,mcg,local}) and (\ref{eq:fem,mdg,local}) using quadrature.

\subsection{Adaptivity}

To achieve the goals stated at the beginning of this section,
the adaptive algorithm chooses individual time steps for the different components based
on the a posteriori error estimates.
Using for example a standard \emph{PID} regulator from
control theory, we choose the individual time steps for each component to satisfy
\begin{equation}
	S_i C k_{ij}^{p_{ij}}r_{ij} = \mathrm{TOL} / N,
	\label{eq:regler,1}
\end{equation}
or, taking the logarithm with $C_i = \log(\mathrm{TOL}/(N S_i C))$,
\begin{equation}
	p_{ij} \log k_{ij} + \log r_{ij} = C_i,
	\label{eq:regler,2}
\end{equation}
with maximal time steps $\{k_{ij}\}$, following
work by S\"oderlind and coworkers \cite{soderlind:pid,soderlind:control}.
Here, $p_{ij}=q_{ij}$ for the \mcgq\ method and $p_{ij}=q_{ij}+1$ for the \mdgq\ method.

To solve the dual problem (\ref{eq:dual,continuous}), which is
needed to compute the stability factors, it would seem that we need to
know the error $e(T)$, since this is used as an initial value for the
dual problem. However, we know from experience that the stability
factors are quite insensitive to the choice of initial data for the
dual problem. (A motivation of this for parabolic problems is given in
\cite{logg:submit:parab}.) Thus in practice, a random (and normalised) value is
chosen as initial data for the dual problem. Another approach is to
take the initial value for the dual problem to be
$\phi(T)=(0,\ldots,0,1,0,\ldots,0)$, i.e., a vector of zeros
except for a single component which is of size one. This gives an
estimate for a chosen component of the error. By other choices of data
for the dual problem, other functionals of the error can be
controlled. In either case, the stability factors are computed
using quadrature from the computed dual solution.

The adaptive algorithm can be expressed as follows: Given a tolerance $\mathrm{TOL}>0$,
make a preliminary estimate for the stability factors and then
\begin{enumerate}
\item
  Solve the primal problem with time steps based on (\ref{eq:regler,1});
\item
  Solve the dual problem and compute the stability factors;
\item
  Compute an error bound $E$ based on (\ref{eq:representation,cg}) or (\ref{eq:representation,dg});
\item
  If $E \leq \mathrm{TOL}$ then stop; if not go back to (i).
\end{enumerate}
Note that we use the error representations
(\ref{eq:representation,cg}) and (\ref{eq:representation,dg})
to obtain sharp error estimates. On the other hand, the error estimates
(\ref{eq:estimate,cg}) and (\ref{eq:estimate,dg}) are used to
determine the adaptive time step sequences.

To limit the computational work, it is desirable that only a few
iterations in the adaptive algorithm are needed. In the simplest case, the error estimate will
stay below the given tolerance on the first attempt. Otherwise, the
algorithm will try to get below the tolerance the second time. It is
also possible to limit the number of times the dual problem is solved.
It should also be noted that to obtain an error estimate at a time
$t=\bar{t}$, different from the final time $T$, the dual problem has
to be solved backwards also from time $\bar{t}$. This may be necessary
in some cases if the stability factors do not grow monotonically as
functions of the final time $T$, but for many problems the stability
factors grow with $T$, indicating accumulation of errors.

Our experience is that automatic computation based on this adaptive
strategy is both reliable (the error estimates are quite close to the
actual error) and efficient (the additional cost for solving the dual
problem is quite small). See \cite{logg:lic:II} for a discussion on
this topic.

\section{Iterative Methods for the Nonlinear System}
\label{sec:iterative}

The nonlinear discrete algebraic equations given by
the \mcgq\ and \mdgq\ methods presented in Section \ref{sec:methods} (including numerical
quadrature) to be solved on every local interval $I_{ij}$ take the
form
\begin{equation}
   \xi_{ijm} =
   \xi_{ij0} +
   k_{ij} \sum_{n=0}^{q_{ij}} w_{mn}^{[q_{ij}]} \ f_i(U(\tau_{ij}^{-1}(s_{n}^{[q_{ij}]})),\tau_{ij}^{-1}(s_n^{[q_{ij}]})),
   \label{eq:equations}
\end{equation}
for $m = 0,\ldots,q_{ij}$,
where $\{\xi_{ijm}\}_{m=0}^{q_{ij}}$ are the degrees of freedom to be determined for component $U_i(t)$ on the interval $I_{ij}$,
$\{w_{mn}^{[q_{ij}]}\}_{m=0,n=0}^{q_{ij}}$ are weights,
$\tau_{ij}$ maps $I_{ij}$ to $(0,1]$: $\tau_{ij}(t) = (t-t_{i,j-1})/(t_{ij}-t_{i,j-1})$, and
$\{s_n^{[q_{ij}]}\}_{n=0}^{q_{ij}}$ are quadrature points defined on $[0,1]$.

\begin{figure}[htbp]
	\begin{center}
		\leavevmode
		\psfrag{0}{$T_{n-1}$}
		\psfrag{T}{$T_n$}
    	\includegraphics[width=12cm,height=6cm]{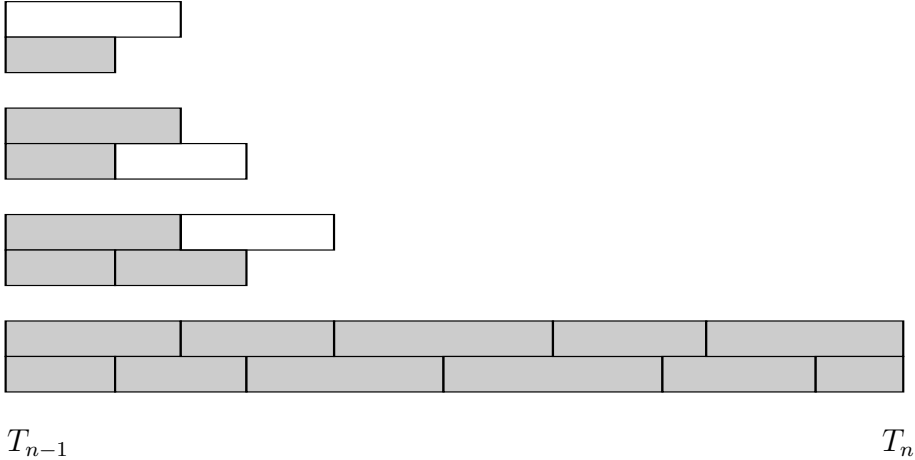}
    	\caption{Multi-adaptive time-stepping within a time slab for a system with two components.}
    	\label{fig:stepping}
  	\end{center}
\end{figure}

The strategy we use to solve the discrete equations
(\ref{eq:equations}) is by direct fixed point iteration, possibly in
combination with a simplified Newton's method. To evolve the system,
we need to collect the degrees of freedom for different components
between two time levels and solve the discrete equations for these
degrees of freedom. We refer to such a collection of elements between
two time levels as a time slab (see Figure \ref{fig:notation}).
New time slabs are formed as we evolve the system starting at time
$t=0$, in the same way as new time intervals are formed in a standard
solver which uses the same time steps for all components.
On each time slab, we thus compute the degrees of freedom
$\{\xi_{ijm}\}_{m=0}^{q_{ij}}$ for each element within the time slab
using (\ref{eq:equations}), and repeat the iterations until the
computational error is below a given tolerance for the computational
error. The iterations are carried out in order, starting at the
element closest to time $t=0$ and continuing until we reach the last
element within the time slab. This is illustrated in Figure
\ref{fig:stepping}.

The motivation for using direct fixed point iteration, rather than
using a full Newton's method, is that we want to avoid forming the
Jacobian (which may be very large, since the nonlinear system to be
solved is for the entire time slab) and also avoid solving the
linearised system. Instead, using the strategy for adaptive damping of
the fixed point iterations as described in the next section, the
linear algebra is built into the adaptive solver. Since often only a
small number of fixed point iterations are needed, typically only two
or three iterations, we believe this to be an efficient approach.

\section{Stiff Problems}
\label{sec:stiff}

As discussed in the previous section, the nonlinear discrete equations
given by the (implicit) multi-adaptive Galerkin methods are solved
using fixed point iteration on each
time slab. For stiff problems these iterations may fail to converge.
We now discuss a simple way to stabilise a stiff system, in order to
make the explicit fixed point iterations convergent.

For simplicity, we assume that the time step sequence, $k_1, k_2,
\ldots, k_M$, is the same for all components.

\subsection{The test equation}
\label{sec:testequation}

To demonstrate the main idea, we consider the stabilisation of the explicit
Euler method applied to the simple {\em test equation}:
\begin{equation}\label{testequation}
	\left\{
	\begin{array}{rcl}
		\dot u(t)+\lambda u(t)&=&0\quad\mbox{for }t>0,\\
		u(0)&=&u_0,
	\end{array}
	\right.
\end{equation}
where $\lambda> 0$ and $u_0$ is a given initial condition. The solution is given by $u(t)=\exp(-\lambda t)u_0$.

The explicit Euler method for the test equation reads
\[
U^n=U^{n-1}-k_n\lambda U^{n-1}=(1-k_n\lambda )U^{n-1}.
\]
This method is conditionally stable, with stability guaranteed if $k_n\lambda\le 2$.
If $\lambda$ is large, this is too restrictive outside transients.

Now, let $K$ be a large time step satisfying $K\lambda >2$ and let $k$ a small time step chosen so that
$k\lambda <2$. Consider the method
\begin{equation}\label{explimplE0}
U^n=(1-k\lambda )^{m}(1-K\lambda )U^{n-1},
\end{equation}
corresponding to one explicit Euler step with large time step $K$ and $m$
explicit Euler steps with small time steps $k$, where $m$ is a positive
integer to be determined. Altogether this corresponds to a time step
of size $k_n=K+mk$. For the overall method to be stable, we require that
$|1-k\lambda|^m (K\lambda-1)\le 1$, that is
\begin{equation}
	m\ge \frac{\log(K\lambda -1)}{-\log|1-k\lambda|}\approx \frac{\log(K\lambda )}{c},
	\label{eq:m}
\end{equation}
if $K\lambda \gg 1$ and $c=k\lambda$ is of moderate size, say $c=1/2$.

We conclude that $m$ will be quite small and hence the small time steps will be used
only in a small fraction of the total time interval, giving a large effective time step.
To see this, define the \emph{cost} as
$\alpha = \frac{1+m}{K+km} \in (1/K,1/k)$,
i.e., the number of time steps per unit interval.
Classical stability analysis gives $\alpha = 1/k = \lambda/2$
with a maximum time step $k=2/\lambda$.
Using (\ref{eq:m}) we instead find
\begin{equation}
	\alpha \approx \frac{1 + \log(K\lambda)/c}{K + \log(K\lambda)/\lambda} \approx
	\frac{\lambda}{c} \log(K\lambda)/(K\lambda) \ll \lambda/c,
\end{equation}
for $K\lambda \gg 1$.
The cost is thus decreased by the cost reduction factor
\[
	\frac{2\log(K\lambda)}{cK\lambda}
\sim \frac{\log(K\lambda)}{K\lambda},
\]
which can be quite significant for large values of $K\lambda$.

\subsection{The general non-linear problem}

For the general nonlinear problem (\ref{eq:u'=f}), the
gain is determined by the distribution of the eigenvalues of the
Jacobian, see \cite{logg:article:stiffode}. The method of stabilising the system using a couple
of small stabilising time steps is best suited for systems
with a clear separation of the eigenvalues into small and large eigenvalues,
but even for the semi-discretised heat equation (for which we have a whole
range of eigenvalues) the gain can be substantial, as we shall see below.

\subsection{An adaptive algorithm}
\label{sec:algorithm}

In \cite{logg:article:stiffode} we present an adaptive algorithm in which both the size of the small
stabilising time steps and the number of such small time steps are automatically determined.
Using adaptive stabilisation, the damping is targeted precisely at the current unstable
eigenmode, which as a consequence allows efficient integration also of problems with
no clear separation of its eigenvalues.

\section{Numerical Examples}
\label{sec:examples}

The numerical examples presented in this section are divided into two categories:
examples illustrating the concept of multi-adaptivity and examples illustrating
explicit time-stepping (or explicit fixed point iteration) for stiff problems.

\subsection{Multi-adaptivity}

The two examples presented below are taken from \cite{logg:article:ma2}, in which
further examples are presented and discussed in more detail.

\subsubsection{A mechanical multi-scale system}
\label{sec:testproblem}

\begin{figure}[H]
	\begin{center}
		\includegraphics[width=12cm,height=2cm]{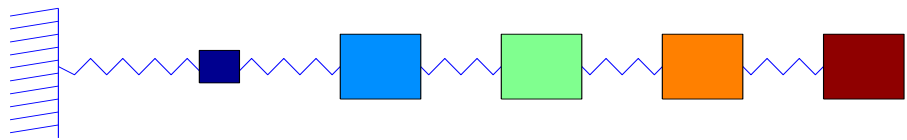}
		\caption{A mechanical system consisting of $N=5$ masses attached together
					with springs.}
		\label{fig:multiadaptivity-system}
	\end{center}
\end{figure}

\begin{figure}[htbp]
        \begin{center}
                \psfrag{e}{\footnotesize \hspace{-0.8cm}$\|e(T)\|_{\infty}$}
                \psfrag{cpu time}{\footnotesize \hspace{-1.2cm}\textsf{cpu time / seconds}}
                \psfrag{function evaluations}{\footnotesize \hspace{-0.7cm}\textsf{function evaluations}}
                \psfrag{steps}{\small \hspace{-0.2cm} \textsf{steps}}
                \includegraphics[width=12cm,height=7.0cm]{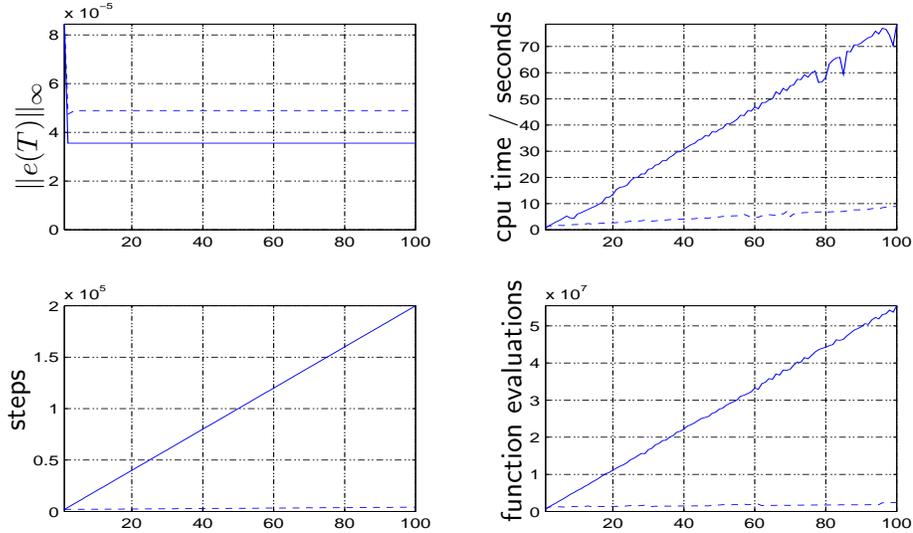}
                \caption{Error, cpu time, total number of steps, and number of function evaluations as
                                        function of the number of masses, for the multi-adaptive $\mathrm{cG}(1)$ method (dashed)
                                        and the standard $\mathrm{cG}(1)$ method (solid).}
                \label{fig:multiadaptivity-results}
        \end{center}
\end{figure}

To demonstrate the potential of the multi-adaptive methods, we consider
a dynamical system in which a small part of the system
oscillates rapidly.
The problem is to compute accurately the positions (and velocities) of
the $N$ point-masses attached together with springs of equal stiffness as in Figure
\ref{fig:multiadaptivity-system}.
We choose a small time step for the smallest mass and large time steps for the
larger masses, and measure the work for the $\mathrm{mcG}(1)$ method as we increase the number of larger masses.
The work is then compared to the work required for the standard $\mathrm{cG}(1)$
method using the same (small) time step for all masses.
As is evident in Figure \ref{fig:multiadaptivity-results}, the work (in terms of function evaluations)
increases linearly for the standard method, whereas for the multi-adaptive method it remains practically
constant.

\subsubsection{Reaction--diffusion}

Next consider the following system of PDEs:
\begin{equation}
	\left\{
		\begin{array}{rcl}
			\dot{u}_1 - \epsilon u_1'' &=& -u_1 u_2^2, \\
			\dot{u}_2 - \epsilon u_2'' &=& u_1 u_2^2,
		\end{array}
	\right.
\end{equation}
on $(0,1)\times(0,T]$ with $\epsilon=0.001$, $T=100$ and
homogeneous Neumann boundary conditions at $x=0$ and $x=1$, which
models isothermal auto-catalytic reactions (see \cite{robban:phd}): $A_1 + 2A_2 \rightarrow A_2 + 2A_2$.
As initial conditions, we take $u_1(x,0) = 0$ for $0 < x < x_0$ , $u_1(x,0) = 1$ for $x_0 \leq x < 1$, and
$u_2(x,0)=1-u_1(x,0)$ with $x_0=0.2$. An initial reaction where substance $A_1$
is consumed and substance $A_2$ is formed will then take place at
$x=x_0$, resulting in a decrease in the concentration $u_1$ and an
increase in the concentration $u_2$. The reaction then propagates to the right until all of substance $A_1$ is
consumed and we have $u_1=0$ and $u_2=1$ in the entire domain.

\begin{figure}[htbp]
	\begin{center}
		\leavevmode
		\psfrag{x}{$x$}
		\psfrag{U1}{$U_1$}
		\psfrag{U2}{$U_2$}
		\psfrag{k1}{$k_1$}
		\psfrag{k2}{$k_2$}
		\includegraphics[width=12cm,height=7.0cm]{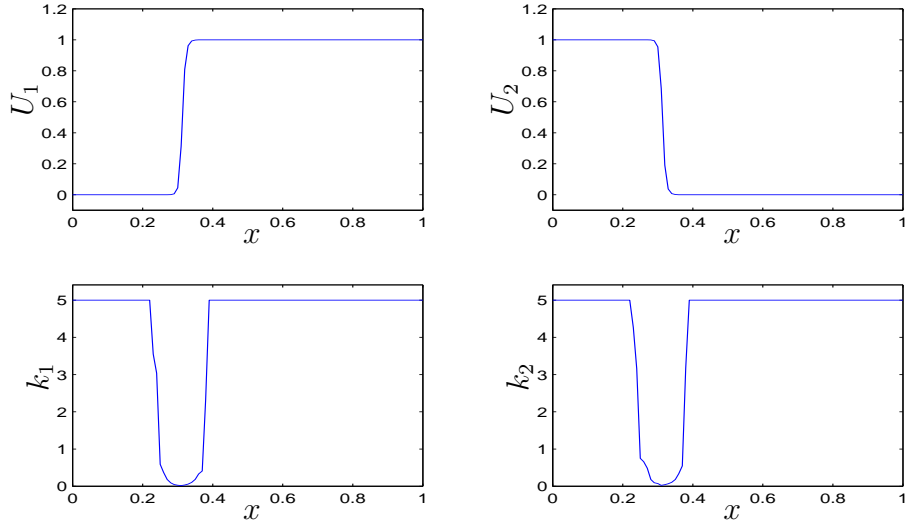}
		\caption{The concentrations of the two species, $U_1$ and $U_2$, at time $t=50$ as function of space (above),
					and the corresponding time steps (below).}
		\label{fig:front-solution}
	\end{center}
\end{figure}

Computing the solution using the $\mathrm{mcG}(2)$ method, we find that the time steps are automatically chosen to
be small only in the vicinity of the reaction front, see Figure \ref{fig:front-solution}, and during the computation
the region of small time steps will propagate to the right at the same speed as the reaction front.

\subsection{Explicit time-stepping for stiff problems}

To illustrate the technique of stabilisation for stiff problems, we present
below some examples taken from \cite{logg:article:stiffode}. In these examples, the cost $\alpha$ is compared
to the cost $\alpha_0$ of a standard implementation of the $\mathrm{cG}(1)$ method in which
we are forced to take a small time step all the time. (These small time steps are marked by
dashed lines in the figures.) Comparison has not been made with an
implicit method, since it would be difficult to make such a comparison
fair; one could always argue about the choice of linear solver and
preconditioner. However, judging by the modest restriction of the average
time step size and the low cost of the explicit method, we believe our
approach to be competitive also with implicit methods, although this
remains to be seen.

\subsubsection{The test equation}

\begin{figure}[htbp]
	\psfrag{t}{$t$}
	\psfrag{U}{$U(t)$}
	\psfrag{k}{$k(t)$}
	\psfrag{5}{}
	\begin{center}
		\includegraphics[width=12cm,height=7.0cm]{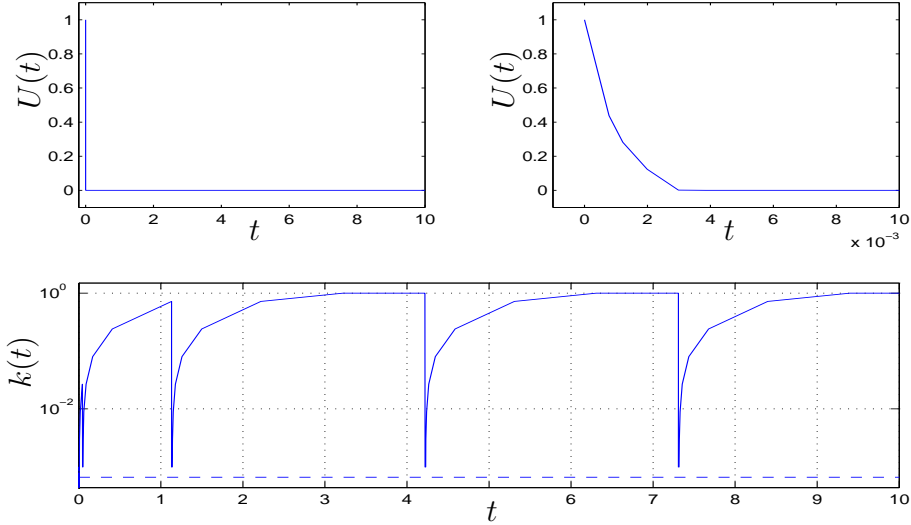}
		\caption{Solution and time step sequence for eq. (\ref{eq:problem1}),
		$\alpha/\alpha_0 \approx 1/310$.}
		\label{fig:problem1}
	\end{center}
\end{figure}

The first problem we try is the test equation:
\begin{equation}
	\label{eq:problem1}
	\left\{
	\begin{array}{rcl}
		\dot u(t)+\lambda u(t)&=&0\quad\mbox{for }t>0,\\
		u(0)&=&u_0,
	\end{array}
	\right.
\end{equation}
on $[0,10]$, where we choose $u_0=1$ and $\lambda=1000$. As is shown in
Figure \ref{fig:problem1}, the time step is repeatedly decreased to stabilise the stiff system,
but overall the effective time step is large and the cost reduction factor is $\alpha/\alpha_0\approx 1/310$.

\subsubsection{The test system}

\begin{figure}[htbp]
	\psfrag{t}{$t$}
	\psfrag{U}{$U(t)$}
	\psfrag{k}{$k(t)$}
	\psfrag{5}{}
	\psfrag{0.05}{}
	\begin{center}
		\includegraphics[width=12cm,height=7.0cm]{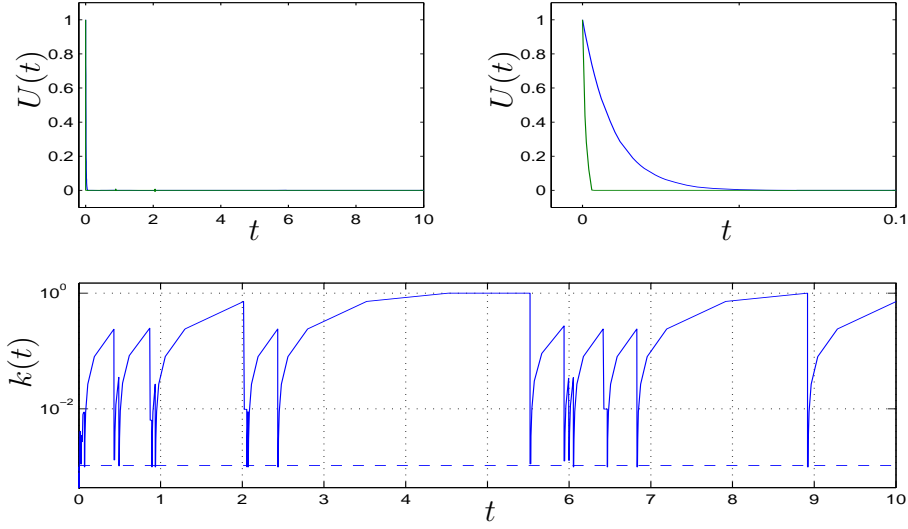}
		\caption{Solution and time step sequence for eq. (\ref{eq:problem2}),
					$\alpha/\alpha_0 \approx 1/104$.}
		\label{fig:problem2}
	\end{center}
\end{figure}

For the test system,
\begin{equation}
	\label{eq:problem2}
	\left\{
	\begin{array}{rcl}
		\dot u(t)+Au(t)&=&0\quad\mbox{for }t>0,\\
		u(0)&=&u_0,
	\end{array}
	\right.
\end{equation}
on $[0,10]$, we take $A = \diag(100,1000)$ and $u_0=(1,1)$.
As seen in Figure \ref{fig:problem2}, most of
the stabilising steps are chosen to damp out the eigenmode corresponding to the largest eigenvalue, $\lambda_2=1000$,
but some of the damping steps are targeted at the second eigenvalue, $\lambda_1=100$.
The selective damping is handled automatically by the adaptive algorithm and the cost reduction
factor is again significant: $\alpha/\alpha_0\approx 1/104$.

\subsubsection{The HIRES problem}

\begin{figure}[htbp]
	\psfrag{t}{$t$}
	\psfrag{U}{$U(t)$}
	\psfrag{k}{$k(t)$}
	\begin{center}
	\includegraphics[width=12cm,height=7.0cm]{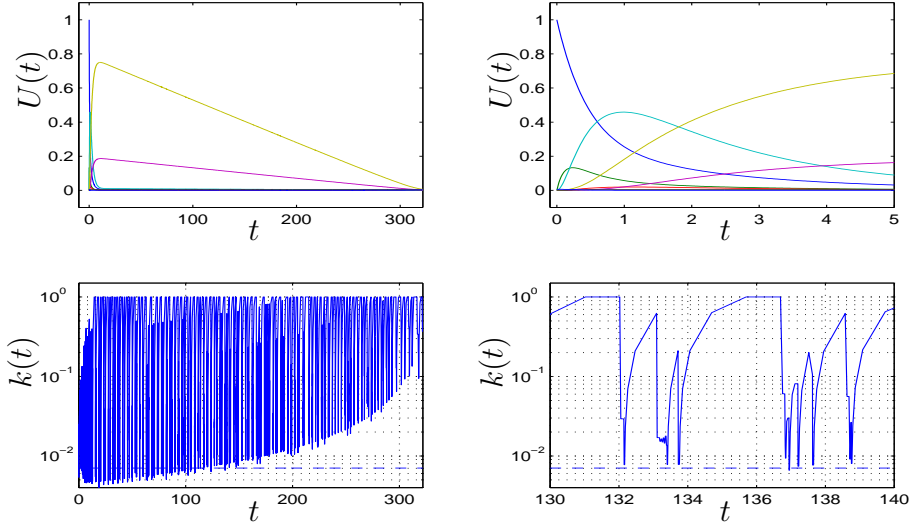}
		\caption{Solution and time step sequence for eq. (\ref{eq:problem6}),
				$\alpha/\alpha_0 \approx 1/33$.}
		\label{fig:problem6}
	\end{center}
\end{figure}

The so-called HIRES problem (``High Irradiance RESponse'') originates from plant physiology and is
taken from the test set of ODE problems compiled by Lioen and de Swart \cite{testset}. The problem
consists of the following eight equations:
\begin{equation}
	\label{eq:problem6}
	\left\{
	\begin{array}{rcl}
		\dot{u}_1 &=& -1.71u_1 + 0.43u_2 + 8.32u_3 + 0.0007, \\
		\dot{u}_2 &=& 1.71u_1 - 8.75u_2, \\
		\dot{u}_3 &=& -10.03u_3 + 0.43u_4 + 0.035u_5, \\
		\dot{u}_4 &=& 8.32u_2 + 1.71u_3 - 1.12u_4, \\
		\dot{u}_5 &=& -1.745u_5 + 0.43u_6 + 0.43u_7, \\
		\dot{u}_6 &=& -280.0u_6u_8 + 0.69u_4 + 1.71u_5 - 0.43u_6 + 0.69u_7, \\
		\dot{u}_7 &=& 280.0u_6u_8 - 1.81u_7, \\
		\dot{u}_8 &=& -280.0u_6u_8 + 1.81u_7,
	\end{array}
	\right.
\end{equation}
on $[0,321.8122]$ (as specified in \cite{testset}). The initial condition is given by
$u_0 = (1.0 , 0 , 0 , 0 , 0 , 0 , 0 , 0.0057 )$.
The cost reduction factor is now $\alpha/\alpha_0 \approx 1/33$, see Figure \ref{fig:problem6}.

\subsection{The heat equation}

\begin{figure}[htbp]
	\psfrag{t}{$t$}
	\psfrag{U}{$U(t)$}
	\psfrag{k}{$k(t)$}
	\psfrag{eigenvalues}{$-\lambda_i$}
	\begin{center}
		\includegraphics[width=12cm]{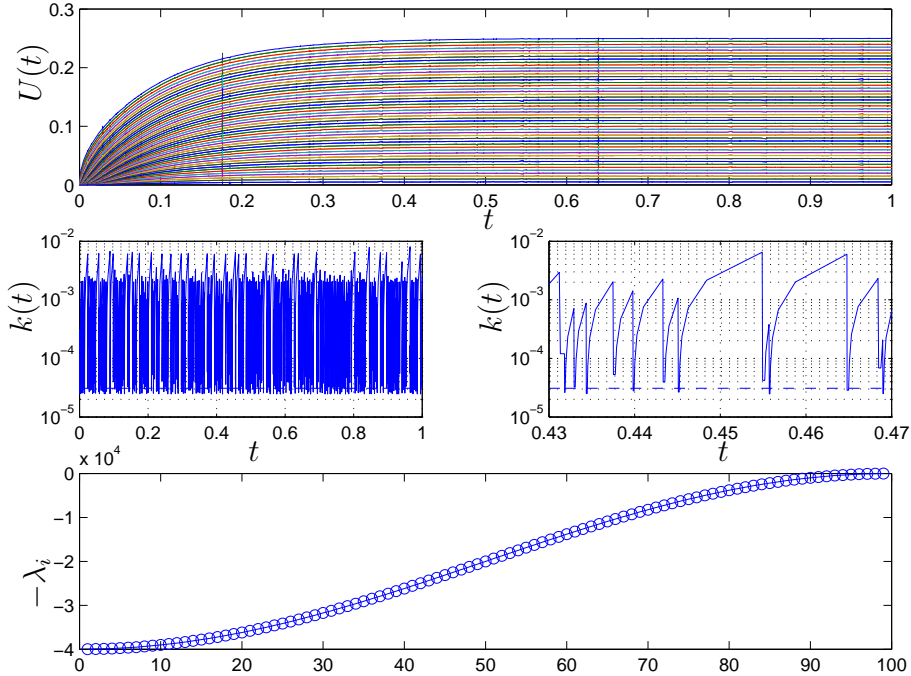}
		\caption{Solution and time step sequence for eq. (\ref{eq:problem10}),
					$\alpha/\alpha_0 \approx 1/17$.}
		\label{fig:problem10}
	\end{center}
\end{figure}

Finally, we consider the heat equation in one dimension:
\begin{equation}
	\left\{
	\begin{array}{rcl}
		\dot{u}(x,t) - u''(x,t) &=& f(x,t), \quad x\in(0,1), \ t>0, \\
		u(0) = u(1) &=& 0, \\
		u(\cdot,t)  &=& 0,
	\end{array}
	\right.
\end{equation}
where we choose $f(x,t)=f(x)$ as an approximation of the Dirac delta
function at $x=0.5$. Discretising in space, we obtain the ODE
\begin{equation}
	\label{eq:problem10}
	\left\{
	\begin{array}{rcl}
		\dot{u}(t) + A u(t) &=& f, \quad t>0,\\
		u(0) &=& 0,
	\end{array}
	\right.
\end{equation}
where $A$ is the \emph{stiffness matrix}. With a spatial resolution of
$h=0.01$, the eigenvalues of $A$ are distributed in the interval
$[0,4\cdot 10^4]$ (see Figure \ref{fig:problem10}). The selective damping
produced by the adaptive algorithm performs well and the cost reduction
factor is $\alpha/\alpha_0\approx 1/17$.

\bibliographystyle{siam}
\bibliography{bibliography}


\newpage
\appendix

\end{document}